# Digital Twins and Their Applications in Modeling Different Levels of Manufacturing Systems: A Review

Sarow Saeedi*

*Department of Industrial, Manufacturing, and Systems Engineering, University of Texas at Arlington, Texas, USA*

## Abstract

Digital Twin (DT) has gained great interest as an innovative technology in Industry 4.0 that enables advanced modeling, simulation, and optimization of service and manufacturing systems. This article provides an extensive review of the literature on digital twins (DTs) and their utilization at the levels of product/production line, production system, and enterprise, and considers how they have been applied under real industrial conditions. This article classifies the types of DTs as well as modeling technologies of DT and applications in different fields, with particular focus on the research of strengths and limitations of Discrete Event Simulation (DES) for systems modelling. A generic structure for DT is proposed, outlining the essential components and flow of data. Case studies demonstrate the benefits of DTs for increased efficiency, reduced downtime, and improved lifecycle management, as well as challenges caused by the complexity of data integration and cybersecurity risk, and high implementation costs. This paper contributes to the growing body of knowledge by identifying both the opportunities and barriers to widespread DT adoption. This study concludes that while DTs offer transformative capabilities for enhancing efficiency and decision-making, overcoming these challenges is crucial for realizing their widespread adoption and impact across service and manufacturing sectors.

---

* Correspondign author. Email address: sarow.saeedi@uta.edu

# 1. Defining the Digital Twin (DT)

In this section, we start by examining various definitions of "Digital Twin" (DT) to identify common themes, and then propose a "best" definition for the scope of this project. Also, we consider how the definition has evolved over the years.

## 1.1. Evolution of the Digital Twin Definition

- **Early Concept (2002-2010):** Although David Gelernter first voiced the digital twin technology through his publication “Mirror Worlds” in 1991 [1], it was Dr. Michael Grieves who first proposed the application of this concept in manufacturing in 2002 and introduced the concept of the digital twin software [2], [3]. Early definitions emphasized the connection between a physical entity and its virtual counterpart for product lifecycle management, primarily focused on design and development stages.

- **Expansion (2010-2015):** As technology advanced, definitions began to incorporate aspects like data analytics, simulation, and real-time connectivity. The exponential increase in computing, storage, and bandwidth capabilities began making the Digital Twin model a reality [3].

- **Current Perspective (2015-Present):** The contemporary view of Digital Twins emphasizes their role in various stages of a product's lifecycle, including operation, maintenance, and end-of-life [4]. Real-time data exchange, simulation, and predictive capabilities are now central to most modern definitions.

## 1.2. Definition from the Literature

The definition of a Digital Twin, which is a term coined by John Vickers of NASA [3], varies across different sources and has evolved. We categorized these definitions into the following categories:

### 1.2.1. General Definitions

- A Digital Twin is a virtual, digital equivalent of a physical product [2].

- DT is a contextualized software model of a real-world object. This means that the behavior of the physical object can be replicated in software and analyzed under the rules that govern its operating environment [5].

- A Digital Twin is a digital representation of a real-world entity or system. This entity can be physical, conceptual, real, or abstract [6].

- It is a virtual representation that serves as the real-time digital counterpart of a physical object or process [7].

- A Digital Twin is a virtual representation of a physical object that evolves and changes over time along with the product it represents [4].

- DT is a digital platform with a role in improving, processing, and managing information at the level of physical and virtual companies [8].

- A Digital Twin is a set of virtual information constructs that fully describe a potential or actual physical manufactured product from the microatomic level to the macro geometrical level [9].

- DT is a set of software components that represents a physical entity, with a data connection between the virtual and physical components, to model the physical entity and to provide services in a given application domain [10].

### 1.2.2. Definitions Focusing on Functionality and Purpose

- DT is a means for **creating a continuum between the physical and virtual worlds**, enabling the transformation of physical objects into programmable entities [5].
- A Digital Twin is a decision-support system; a virtual representation of a physical object [11].
- The Digital Twin is a production system virtual representation that can be executed in various simulation environments, through the synchronization between the real and virtual systems, using mathematical models and appropriate information [12].

### 1.2.3. Industry-Specific Definitions

- In manufacturing, it is a technology that creates a virtual model of a physical product to optimize the production process [7].
- In the smart city sector, a Digital Twin is an urban management technology that collects data for urban problems and provides solutions by reflecting real-world information in a virtual space model that simulates the real space [7].

### 1.2.4. Definitions Focusing on Data and Communications

- The Digital Twin consists of a physical element, a virtual counterpart, and a communication channel (the "Digital Thread") between the two [3].
- Advanced Digital Twins, sometimes referred to as "cognitive" or "intelligent" twins, are empowered by AI and big data and can self-modify and evolve independently of human intervention [13].
- A predictive Digital Twin is characterized by continuous, real-time communication and also includes cybersecurity components and AI algorithms to increase the accuracy of simulations and allow for predictions [14].

## 1.3. Key Concepts

The core themes that appear in almost every definition include:

- **Virtual Representation:** A digital twin is a digital model or representation.
- **Physical Counterpart:** It corresponds to a real-world object, process, or system.
- **Data Connection:** A link between the physical and virtual elements, often involving real-time data.
- **Simulation & Analysis:** The digital model allows simulation, analysis, and prediction.
- **Life Cycle:** The Digital Twin is not just for design; it's useful across all product/process stages.

## 1.4. Working Definition

Based on these findings, the "best" working definition reflecting the scope of this project can be summarized as follows:

*"A Digital Twin is a dynamic, virtual representation of a physical product, process, system, or enterprise, characterized by a bi-directional data connection that enables real-time data exchange between the physical and virtual entities. This representation uses simulation, modeling, and analysis to understand, monitor, predict, optimize, and manage the performance and behavior of the physical counterpart throughout its lifecycle. Furthermore, the Digital Twin, based on continuous data updates, is a dynamic and adaptive model of its physical equivalent."*

## 2. Digital Twins and Industry 4.0

For considering the contribution of DT to Industry 4.0, we used the search key terms of [*"Digital Twin" AND "Industry 4.0"*] in Science Direct, Google Scholar, and Semantic Scholar databases. As Figure 1 shows, 120 papers were found in the first round, and we used the following criteria to screen the papers and select 75 documents in the screening phase:

- **Publisher**: Is the article published in a conference/journal or as a pre-print?
- **Operational Outcomes**: Does the study report specific, quantifiable operational improvements or outcomes from the digital twin implementation?
- **Industry 4.0 Framework**: Does the study explicitly connect to or operate within an Industry 4.0 framework or context?
- **Implementation Evidence**: Does the study include practical implementation evidence beyond purely theoretical frameworks?

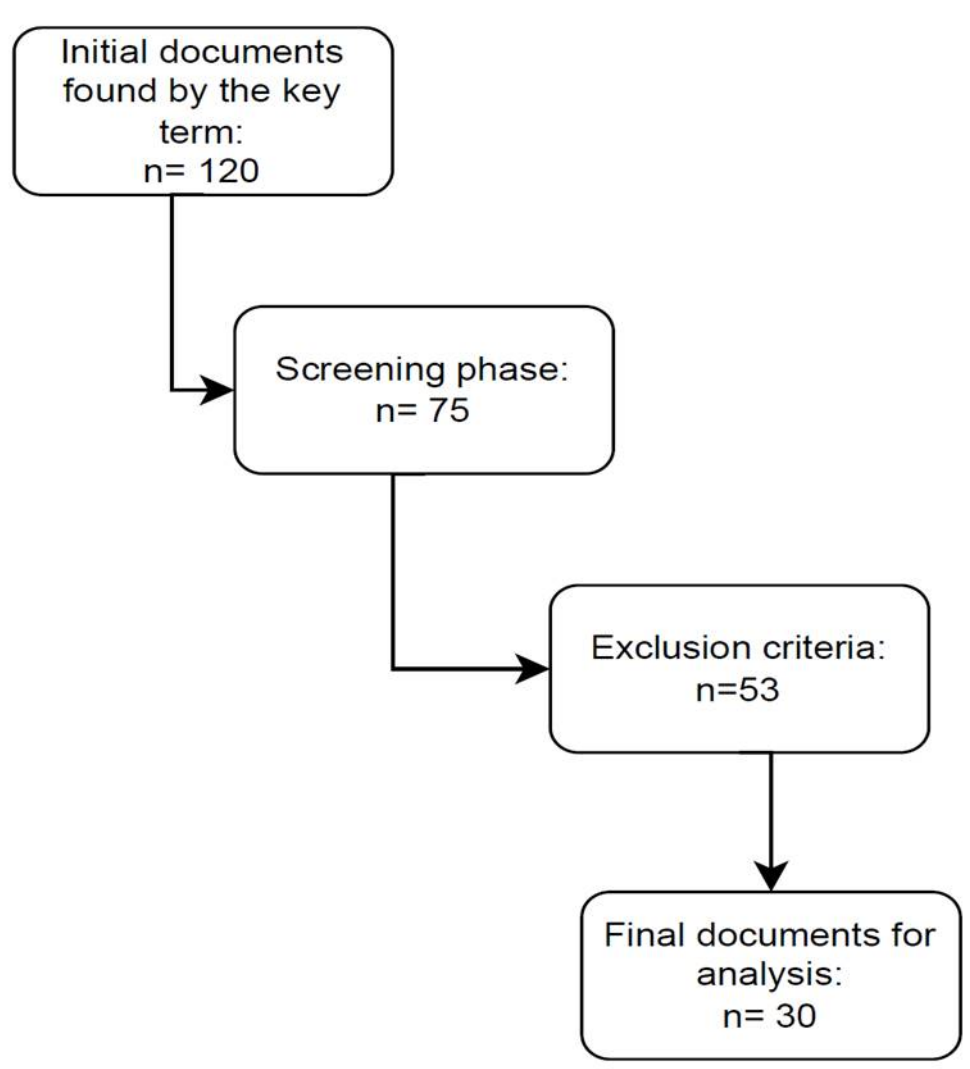

*Figure 1. Systematic review methodology for the contribution of DT to Industry 4.0*

Afterward, the results were scrutinized by reading their abstracts, which excluded some of them for further consideration. The results were refined to those focused merely on the contribution of DT to Industry 4.0. In this phase, 53 documents were selected for further consideration. Then, in the final phase, we utilized the following criteria to select the final documents:

- **Digital Twin Type and Purpose**: Identifying the purpose of the digital twin (e.g., energy optimization, process control, performance measurement), and removing the document if no technological components are described (e.g., software, communication protocols, and modeling approaches).
- **Digital Twin Implementation Approach**: Prioritizing information in methodology and technical implementation sections to specify whether key technologies are used (e.g., IoT, AI, machine learning, communication protocols).
- **Quantitative Performance Improvements**: Focusing on results and discussion sections, and extract specific quantitative improvements demonstrated by the digital twin, and remove the document if no quantitative improvements are explicitly stated

## 2.1. Contribution Levels of DT to Industry 4.0

The literature review shows that the DT contributes to Industry 4.0 across three broad levels. These levels are shown in Figure 2.

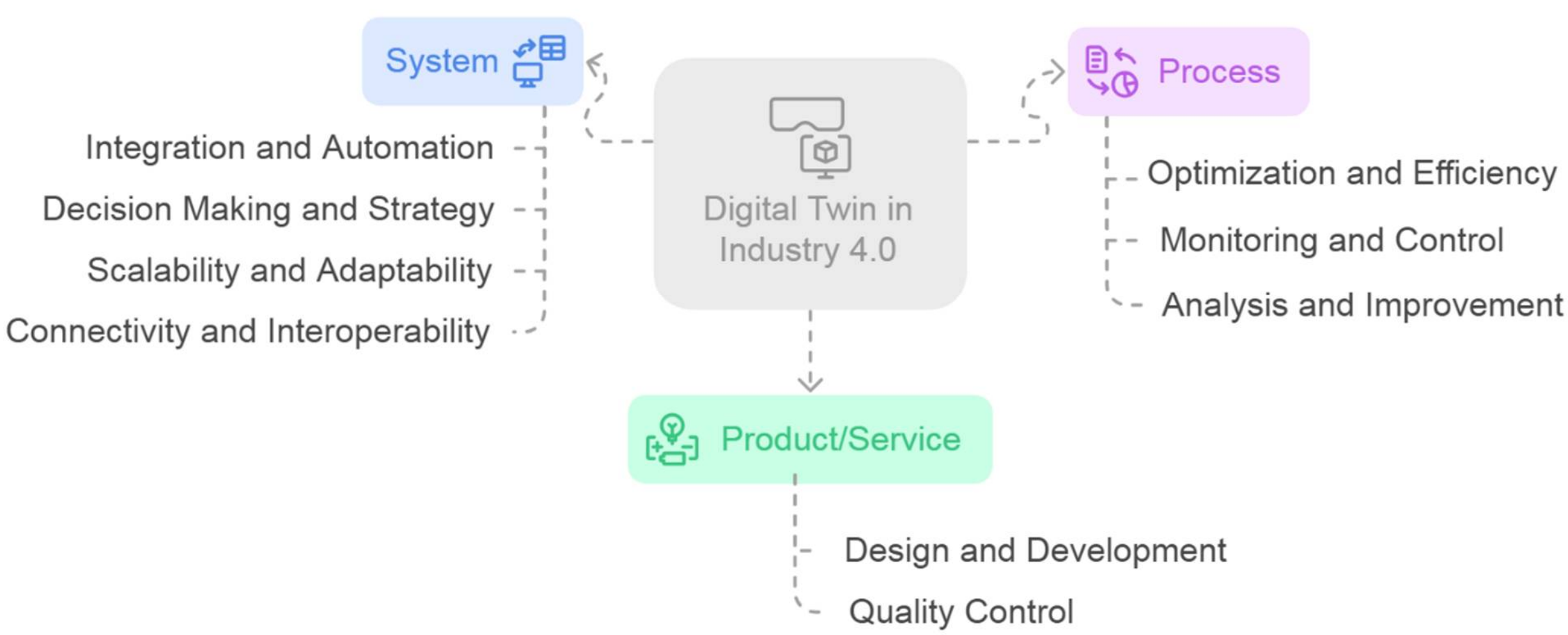

*Figure 2. Levels of Contribution of DT to Industry 4.0 (Graphic generated using Napkin.AI)*

### 2.1.1. Product/Service Level:

- **Design and Development**: *Virtual Prototyping* is an area in which DT can help companies identify and resolve potential issues before manufacturing the final product, reducing engineering costs and improving operational procedures [15], [16], [17]. DT can also enable a more virtual, system-based design approach [16]. A *Comprehensive Digitized Footprint* of the entire product development cycle, from design to deployment, is another area in this category that DT can contribute to Industry 4.0 [16]. Digital Twin also supports *Product Lifecycle Management*, ensuring that products are monitored and optimized throughout their entire lifespan [18]. Finally, Customization and Personalization are another area in this category that is enabled by DT [19].
- **Quality Control**: DT enables *Defect Reduction* by real-time monitoring, analysis, and optimization of production processes [20]. DT can also contribute to Industry 4.0 through *Parts Twinning*, as DT assists engineers and developers in better understanding a given part's mechanical, physical, and intellectual features in the context of the entire product [16].

#### 2.1.2. Process Level:

- **Optimization and Efficiency**: DT enables *Process Optimization* by equipment monitoring, process adjustments, and digital maintenance [15], [16]. It also improves *Resource Management* by assessing material usage, discovering inefficiencies, and optimizing tool-tracking systems. *Real-time adaptability* is another area in this category in which DT contributes to Industry 4.0 by creating virtual representations that adapt to changes in the physical environment [16].

- **Monitoring and Control**: DT allows for *Remote Monitoring* and operation of equipment and systems [16]. It also helps *Anomaly Detection* by preventing and controlling possible failures encountered during the use of real objects [17]. Furthermore, DT facilitates *Predictive Maintenance* to minimize downtime and extend the lifespan of critical equipment [21].

- **Analysis and Improvement**: *Bottleneck Detection* is a major contribution of DT to Industry 4.0 in this category. Digital Twin detects, diagnoses, and improves bottleneck resources using utilization-based bottleneck analysis, process mining, and diagnostic analytics [22], [23], [24]. Also, DT offers manufacturers valuable *Data-Driven Insights* that were previously unattainable, driving innovation and operational efficiency [25].

#### 2.1.3. System Level:

- **Integration and Automation**: Digital Twin enables the *Seamless Integration* of distinct production processes, from planning to actuators [16]. It can also support smart *Autonomous Systems* that monitor and control different machinery, tools, robots, and automated vehicles [18]. Furthermore, DT encompasses *Automation*, data interchange, and manufacturing processes [16].

- **Decision-Making and Strategy**: On the one hand, DTs can lead to higher performance in Industry 4.0 as they can make *Informed Judgements*. On the other hand, by improving the decision-making resulting from wise judgments, DTs can enhance *Strategic Planning* at the organizational level [16], [26], [27].

- **Scalability and Adaptability**: Due to their *Versatile* nature, DTs have the potential to revolutionize Industry 4.0. So, they must be scalable and adaptable to various processes and products [16], [19].

- **Connectivity and Interoperability**: Digital Twins can integrate with the Industrial Internet of Things (*IIoT*), which leads to the establishment of collaborative frameworks between edge and cloud computing. They can also facilitate *Data Exchange* by connecting physical and digital counterparts [19].

# 3. Types of Digital Twins, Modeling Technologies, & Applications

## 3.1. Types of DTs

According to Singh et. al. [28], there are five general types of DTs, each of which has different subcategories. Figure 3 shows different types of DTs.

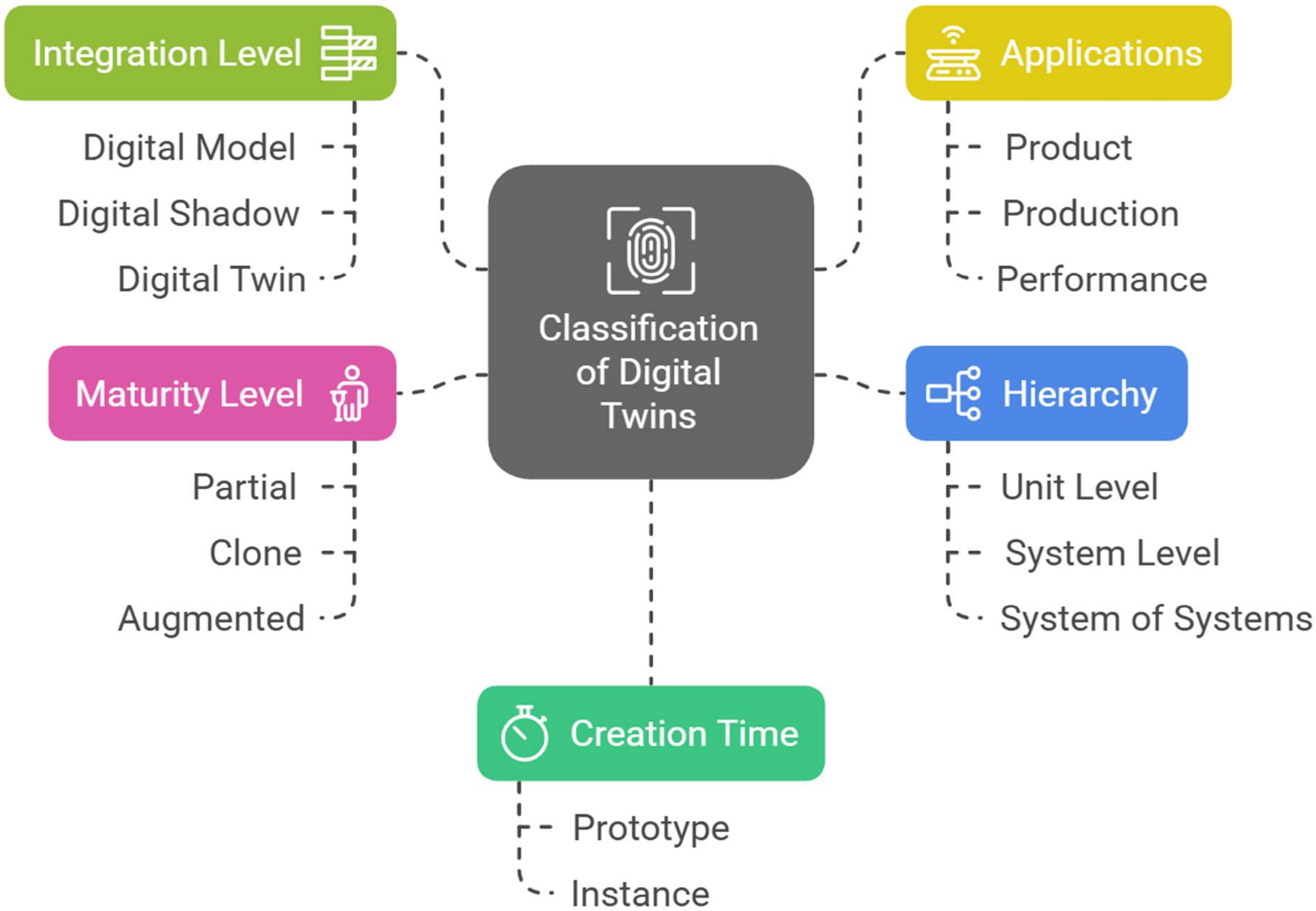

*Figure 3. Types of DTs (Graphic generated using Napkin.AI)*

## 3.2. Modeling Technologies to Construct DTs

As mentioned by Tao et. al. [29], there are four types of models for constructing a digital twin. These models are shown in Figure 4. According to these categories, they introduced a variety of technologies and tools that can be utilized to construct DTs, which are depicted in Figure 5.

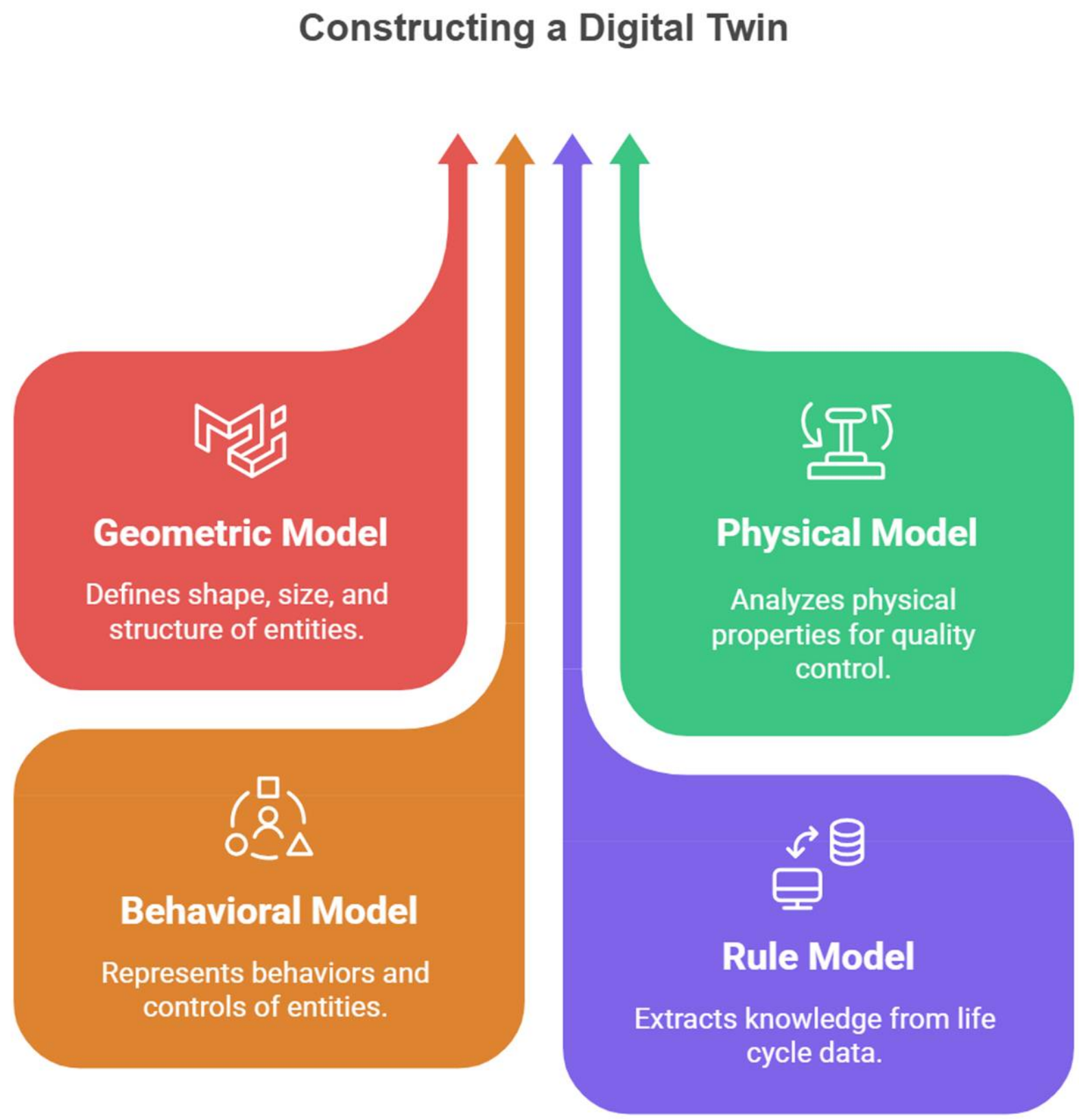

*Figure 4. Types of Models for Constructing DTs (Graphic generated using Napkin.AI)*

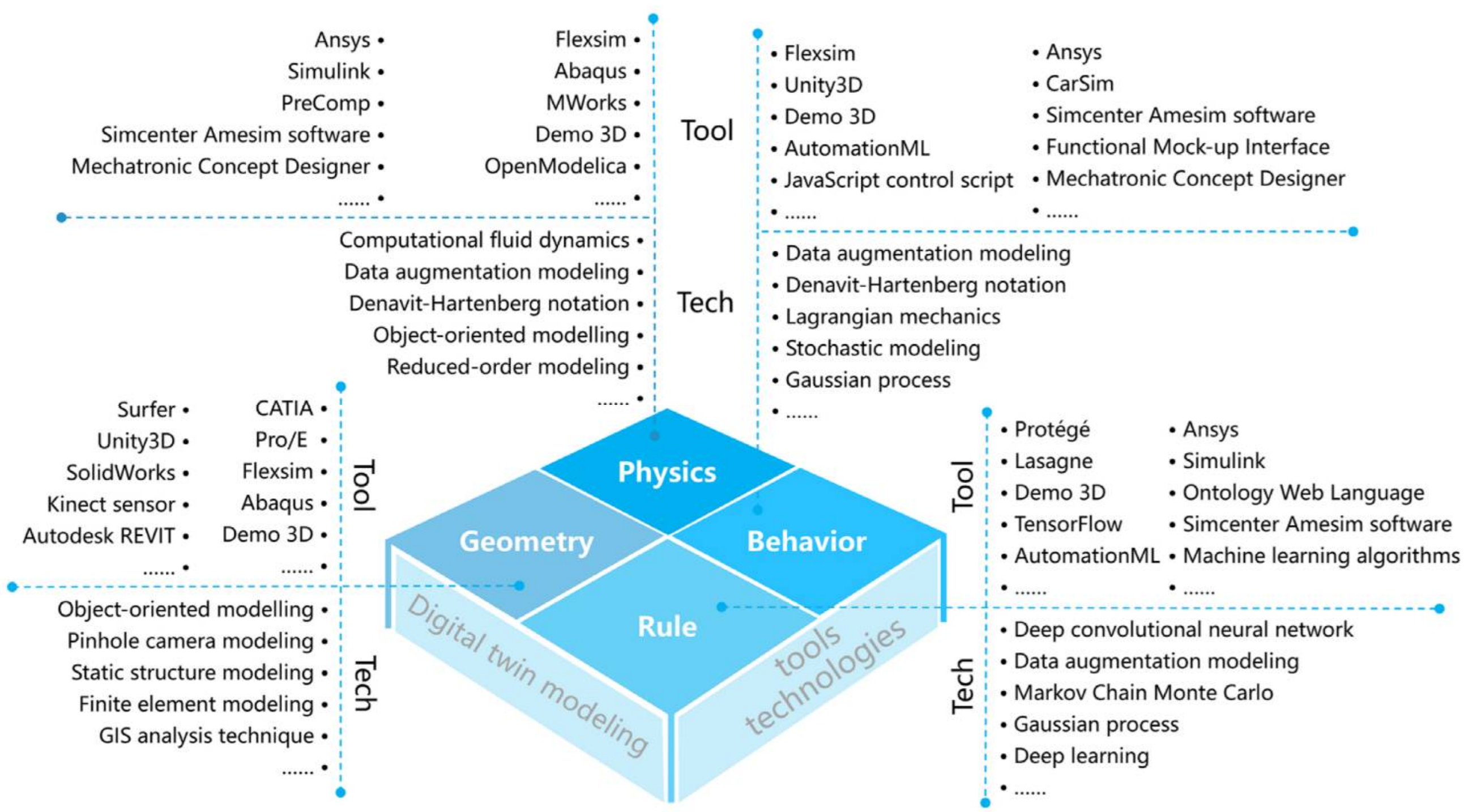

*Figure 5. Tools and Technologies in DTs Model Construction (Adopted from [29])*

## 3.3. Systems Modeled by DTs

Several systems can be modeled with DT. Based on the current literature, Digital Twins are most applied to model *Manufacturing Systems* [25], [29], [30]. Other systems that can be modeled by DTs are *Healthcare Systems*, *Energy Systems*, *Aerospace Systems*, *Transportation and Logistics Systems*, *Engineering Construction Systems*, *Agricultural Systems*, *Software Engineering Systems*, *Infrastructural Systems*, and *Educational Systems*.

DTs have multiple applications in modeling the mentioned systems. These applications are summarized in Table 1.

*Table 1. Applications of DTs across Industries/Systems*

| Industry / System | Application | Description | Resource |
|---|---|---|---|
| **Aerospace and Aeronautics** | Performance and Reliability Optimization | Optimizing the performance and reliability of spacecraft/aircraft to ensure mission success. | [12] [31] |
| | Predictive Maintenance | DTs predict and resolve maintenance issues, reducing downtime and costs. | [12], [31] |
| | Flight Simulation and Scenario Prediction | Flight simulation before launch and future scenarios prediction by continuously mirroring actual flight conditions. | [12], [31] |
| | Damage Diagnosis | Diagnosing damage to spacecraft/aircraft, for in-situ repairs and mission adjustments. | [32] |
| | Parametric Studies and Design Modifications | Studying the effects of modified parameters not considered during design. | [32] |
| | Fan-Blade Reconditioning Automation | Automate reconditioning using vision-based algorithms and robotics for aerospace maintenance. | [33] |
| | Machining Process Mimicking | Mimic of the physical machining processes of aircraft components like rudders. | [34] |
| | Rocket Engine Condition Prediction | Prediction of real flight conditions and their impacts on rocket engine start-up, enhancing reliability. | [35] |
| **Manufacturing** | Product Design and Verification | DTs allow designers to virtually verify product designs, test iterations, and optimize features based on real-time data and consumer feedback. | [36], [37], [38], [39], [40], [41] |
| | Material Selection Optimization | Simulation of manufacturing phases for different materials to optimize material selection based on properties, cost, and environmental impact. | [42] |
| | Production Planning and Control | Digital Twins aid in resource management, production planning, and process control, improving efficiency and decision support. | [36], [38], [40], [41], [43] |
| | Predictive Maintenance and Downtime Reduction | Prediction of failures, enabling scheduled and preventive maintenance to minimize downtime and costs. | [38], [40], [41], [43], [44] |
| | Real-time State Monitoring and Service | Providing real-time monitoring of product operation and facilitating various product services, including energy consumption analysis and user behavior analysis. | [36], [38], [40], [41], [43] |
| | Waste Recovery and Remanufacturing Support | DTs support waste electrical and electronic equipment recovery, aiding manufacturing and remanufacturing operations throughout the product lifecycle. | [45] |
| | Additive Manufacturing (AM) Optimization | Optimization of 3D printing processes, reducing trial-and-error tests, and making AM time- and cost-effective. | [46] |

| Industry / System | Application | Description | Resource |
|---|---|---|---|
| | Defect Detection in AM Parts | Using a combination of in-situ sensing and machine learning. | [47] |
| | Robot Programming and Validation | DTs are used in offline, online, and manual robot programming methods, as well as for validating human-robot collaboration safety standards. | [36], [48], [49] |
| | Bottleneck Detection and Production Optimization | DTs analyze production processes in real-time to identify bottlenecks and inefficiencies, enabling optimization of workflow and resource allocation. | [22], [23], [24], [36], [38], [40], [41], [43] |
| **Healthcare and Medicine** | Radiology Department Optimization | DTs can optimize hospital radiology departments to improve patient care, reduce waiting times, and enhance equipment utilization. | [30] |
| | Personalized/Precision Medicine | DTs facilitate personalized treatments by mimicking patient behavior and suggesting tailored cures based on individual data. | [50], [51], [52] |
| | Organ and Body Part Modeling | Creating digital models of organs and body parts (e.g., brain, heart) for research, disease prediction, and treatment planning. | [30], [52] |
| | Disease Detection and Management | Detecting Ischemic Heart Diseases (IHD) and Alzheimer's by collecting patient data and predicting disease progression. | [30], [52], [53] |
| | Clinical Trial Acceleration | DTs are considered for clinical trials to accelerate medical innovations and regulatory approvals through data-driven predictions. | [30] |
| **Power Generation/Energy** | Renewable Energy Optimization | Through predicting energy output based on weather conditions and optimizing the placement and operation. | [12], [30], [54], [55], [56], [57] |
| | Wind Turbine Performance Improvement | DTs, like WindGEMINI, improve wind turbine performance through predictive maintenance and long-term energy production evaluation. Plus, predicting wind flow patterns in wind farms to optimize turbine placement. | [30], [55], [56] |
| | Optimal Location for Energy Generation | By simulating environmental conditions and energy yield in wind farms or solar plants. | [12], [30], [56] |
| | Nuclear Power Plant (NPP) Improvement | Improving NPP control algorithms, equipment diagnosis, operator training, and lifecycle management. | [12], [30], [56] |
| | Issue Identification and Prediction | Identification of potential failures or inefficiencies, allowing for proactive measures. | [16], [30] |
| | Electricity Distribution Optimization | DTs can simulate different scenarios to improve grid stability, reduce losses, and enhance overall efficiency. | [30] |
| **Automotive** | Predictive Maintenance | DTs predict brake pad wear and facilitate predictive maintenance by comparing real-time and simulated data. | [12], [58] |
| | Personalized Customer Services | DTs enable car manufacturers to provide personalized services based on vehicle operational and behavioral data. | [30] |
| | Vehicle Sales Enhancement | Using AR, DTs provide a 360° view of vehicles and integrate customer preferences. | [30] |
| | Formula 1 Car Performance Analysis | DTs analyze car performance data in Formula 1 racing to improve car design, reliability, safety, and race strategies. | [30] |

| Industry / System | Application | Description | Resource |
|---|---|---|---|
| | Autonomous Vehicle Validation | DTs validate safety and movement algorithms for autonomous vehicles in virtual environments. | [30] |
| **Smart City** | Water Supply Management | Digital Twins manage city water supply efficiently, reducing operational hours and improving service delivery to citizens. | [12], [30], [59] |
| | Urban Planning and Decision Making | Digital Twins aid in urban planning and decision-making, fostering sustainable, economic, and environmentally friendly cities. | [12], [30], [59] |
| | Emergency Response and Disaster Management | DTs are used for disaster management and emergency response planning, predicting floods. | [59], [60], [61], [62] |
| | Infrastructure Analysis and Risk Assessment | Digital Twins allow analysis of city infrastructure in different scenarios and risk assessment. | [30], [59] |
| | Traffic Management and Optimization | DTs manage and optimize traffic flow in smart cities, reducing congestion. | [59] |
| **Transportation and Logistics** | Vessel Performance Monitoring and Diagnostics | DTs monitor vessel performance in real time and diagnose potential failures, improving asset availability and operational readiness. | [12], [30], [31] |
| | Predictive Maintenance | DTs facilitate predictive maintenance for crane vessels. | [30] |
| | Remote monitoring of assets | In this sector, DTs enable remote monitoring of transportation assets, such as vehicles and fleets, to track their location, performance, and maintenance needs. | [16] |
| | Supply chain optimization | By simulating different scenarios to improve efficiency, reduce costs, and enhance delivery times. | [16], [57] |
| | Port Digitalization | Digitalizing ports by gathering real-time environmental and operational data. | [30] |
| **Agriculture** | Livestock Remote Monitoring | Monitoring livestock health, movement, and estrus cycles. | [12], [30] |
| | Pest and Disease Identification | Digital Twins aid in identifying pests and diseases in plants, enabling effective and timely interventions. | [63] |
| | Crop Management Cost | Evaluation of the cost-effectiveness of crop management treatments and tracking machinery in real-time. | [63] |
| | Vertical Farming Implementation | DTs can optimize different parameters in vertical farms to improve productivity. | [64] |
| | Livestock Farm Management | DTs offer frameworks for livestock farm management, including modeling, analysis, simulation, and visualization for optimized operations. | [65] |
| | Sustainability Promotion | Digital Twins track carbon emissions, biodiversity, and water catchment services in agriculture, promoting sustainable practices. | [30] |
| **Education** | Technical Course Enhancement | By representing multiple domains and visualizing system performance, improving understanding, knowledge exchange, and ensuring the safety of students. | [30], [66] |
| | Authentic Learning Experiences | DTs promote authentic learning experiences, enabling effective knowledge construction, skill mastery, and self-efficacy. | [67] |
| | Physical Twin Behavior Learning | DTs allow learning about physical twin behavior in real-world conditions, | [67] |

| Industry / System | Application | Description | Resource |
|---|---|---|---|
| | | providing immediate feedback and problem-solving opportunities. | |
| | Inquiry-Based Learning | Digital Twins facilitate inquiry-based learning during system development and testing. | [67] |
| | Individualized Learning | DTs enable each student to work on an individual DT, offering personalized learning experiences and resource access. | [67] |
| | Distance Learning Support | DTs are great tools for distance learning, providing access to virtual physical twins when physical access is not possible. | [67] |
| **Construction** | Building Information Modeling (BIM) | Enhancing BIM by generating and managing construction project information across its lifecycle. | [12], [57], [68] |
| | Design and Engineering Process Reduction | Digital Twins reduce the overall design process and minimize additional costs during rework by virtually verifying product designs. | [16] |
| | Problem Solving and Data Verification | DTs help in problem-solving by allowing data to be added, modified, and verified against real-life scenarios, improving decision-making. | [30] |
| | Structural System Integrity Assessment | DTs assess structural system integrity and ensure buildings do not fail under applied forces, promoting modular construction. | [30] |
| | Project Management and Scheduling | Digital Twins assist in managing resources, materials, schedules, and quality during construction, optimizing project execution. | [30] |
| | Building Asset Condition Monitoring | Digital Twins continuously monitor building asset conditions, optimizing maintenance and services, and aiding predictive maintenance. | [57] |
| | Demolition and Heritage Conservation | Digital Twins help in problem-solving for future projects, conserve heritage assets virtually, and identify potential hazards during demolition. | [30] |
| **Retail** | Tailored Customer Experience | Providing suggestions based on customer interest patterns, enhancing satisfaction. | [30] |
| | Logistics and Supply Chain Optimization | DTs optimize retail logistics and supply chains, improving inventory planning, demand forecasting, and real-time monitoring. | [12], [69] |
| | Store Layout Planning | DTs facilitate testing the efficacy of different store layouts before implementation, optimizing store design for better customer flow and sales. | [70] |
| | Emergency and Disruption Decision-Making | Digital Twins simulate "what-if" scenarios during emergencies or disruptions like COVID-19, aiding in decision-making for retailers. | [71] |

### 3.3.1. The Best DT model at the Product, Production Facility, and Enterprise Level:

The most suitable methods and tools vary depending on the specific context, including factors like the industry type (manufacturing or service), system complexity, and compatibility with existing infrastructure (such as ERP[1] or MES[2] systems). For example, a company utilizing Siemens ecosystems may find MindSphere ideal for seamless integration at all levels, while a company prioritizing cloud scalability might choose Azure Digital Twins at the enterprise level. Nevertheless, the recommendations provided in this section are considered versatile options for analyzing Service and Manufacturing Systems using Digital Twins.

*DT Models, Tools, and Technologies at the Product/Production Level:*

DT plays various roles throughout the product lifecycle, including design, manufacturing, delivery, usage, and end-of-life [72]. The primary goal at this level is to model and simulate the behavior of individual components or processes to optimize design, performance, and efficiency.

Researchers have utilized all models mentioned in Figure 4 to construct a product/production line DT (e.g., [73], [74], [75]). The most mentioned models and technologies in the literature for DT at the Product/Production level include *3D Modeling and Visualization* (e.g., [72], [74]), *Finite Element Modeling (FEM)* (e.g., [76]), *Computational Fluid Dynamics (CFD)* (e.g., [77]), *Big Data Analytics (BDA)*, *Machine Learning Models* (e.g., [72]), and *Discrete Event Simulation (DES)* (e.g., [73] for constructing both product and process DT). To construct a product DT, different tools can be utilized. For instance, *Siemens NX*, *AutoDesk Fusion 360*, *Ansys*, *SolidWorks*, and *CATIA* for creating the Virtual Replica of the product DT [72], [78] and production line DT [75], *MATLAB*[3], *Simul8*[4], and *FlexSim*[5] for the simulation environment of the product DT, *Automation Markup Language (AutomationML)* [79], *DDS, HLA/RTI, and MQTT* for the data transfer and communication under the data component of the product DT [72], and finally *Siemens Tecnomatix* for plant simulation [24], [73].

*DT Models, Tools, and Technologies at the Production Facility Level:*

This level involves modeling entire production facilities, encompassing multiple production lines and systems. The aim is to optimize the overall facility performance, including energy management, layout optimization, and resource allocation.

All models mentioned in Figure 4 to construct a DT have been utilized by researchers. Most researchers focused on plant and shop-floor DTs (e.g., [75], [80], [73], [74]). Models and technologies utilized at this level consist of *Discrete Event Simulation (DES)* [73], *Finite Element Modeling (FEM)* [75], *Neural Networks (NN)* [75], *Object-oriented Modeling* [80], *Machine Learning Models* [73], [81], and *Big Data Analytics* [81], [82]. Also, the most mentioned tools for constructing a DT at the production facility level are as follows: *3D Max*, *CAX*[6] [75], *CATIA* [75], *MES* [83], *ERP* [75], *CAXperts UniversalPlantViewer* [84], and *Siemens Tecnomatix* [73]. Also, *OPC Unified Architecture (OPC-UA)* and *MTConnect* protocols for communication [73], [74], and *Siemens Mindsphere*, *Thingworxs Kepware*, and *Vuforia Spatial Toolbox* for integration of different technologies [74].

*DT Models, Tools, and Technologies at the Enterprise Level:*

At the enterprise level, Digital Twins are used to model the entire organization, including multiple facilities, supply chains, and business processes. The objective is to achieve visibility and optimize strategic decision-making [16]. Referring to Figure 3, the DTs that are utilized for modeling at the enterprise level fall under the System of Systems (SoS) level from the Hierarchy classification, meaning that the DT at the SoS level integrates various stages of the product across its entire life cycle [28]. Few researchers have considered modeling DTs at this level [73].

---

[1] Enterprise Resource Planning
[2] Manufacturing Execution System
[3] https://www.mathworks.com (Digital Twin | What Is a Digital Twin? – Accessed: March 2025)
[4] https://www.simul8.com (Build rapid, cost-effective digital twins with Simul8 – Accessed: March 2025)
[5] https://www.flexsim.com (FlexSim + Digital Twin – Accessed: March 2025)
[6] Computer-Aided eXpert

The most utilized model by the researchers to construct a DT model at the enterprise level is *Reference Architectural Model I4.0 (RAM I4.0)* [85], which is a model proposed in the German government's Industry 4.0 program, and it is used by several researchers (e.g., [74], [86], [87]).

Tools that can be employed to build a DT model at the enterprise level are: *ERP*, *Customer Relationship Management (CRM)* [86], *Oracle IoT Digital Twin*, *SAP Digital Twin Framework* [57], and *Siemens Xcelerator* [88]. *AutomationML* [86] is used as a communication tool at the enterprise level.

Table 2 summarizes the recommendations on model, technology, and tool selection for building DTs at different levels.

*Table 2. Recommended Models, Technologies, and Tools for building a DT at different levels*

| **Level** | **Modeling Method / Technology** | **Recommended Tools for Creating the Virtual Replica** | **Recommended Tools for Communication / Integration** |
|---|---|---|---|
| Product/Production Line | 3D Modeling and Visualization; FEM; CFD; BDA; ML Models; DES | Siemens NX; AutoDesk Fusion 360; Ansys, SolidWorks; CATIA; MATLAB; Simul8; FlexSim; Siemens Tecnomatix | AutomationML; DDS; HLA/RTI; MQTT |
| Production Facility | DES; FEM; NN; Object-oriented Modeling; ML Models; BDA | 3D Max; CAX; CATIA; MES; ERP; CAXperts UniversalPlantViewer; Siemens Tecnomatix | OPC-UA; MTConnect protocols; Siemens Mindsphere; Thingworxs Kepware; Vuforia Spatial Toolbox |
| Enterprise | RAM I4.0 | ERP; CRM; Oracle IoT Digital Twin; SAP Digital Twin Framework; Siemens Xcelerator | AutomationML |

## 4. Generic Digital Twin System Architecture

There are several structures and architectures defined for DTs in the literature (e.g., [73], [74], [75], [76], [78], [80], [86], [89], [90], [91]). Based on these structures, we can present a general architecture of DTs. Figure 6 shows this general architecture. In the general DT architecture shown in Figure 6, sensor data from the physical system passes through a data acquisition and network layer into a storage and integration layer, which also incorporates external enterprise data. The aggregated data is then processed by analytics to generate insights that feed into the Digital Twin model, which is built, assembled, fused, verified, and modified iteratively. The outputs from this Digital Twin model (such as simulations and predictions) are displayed on a visualization dashboard and provided to a decision support module, which issues control commands or updates back to the physical system. So, a continuous loop will be formed between the digital replica and the real system.

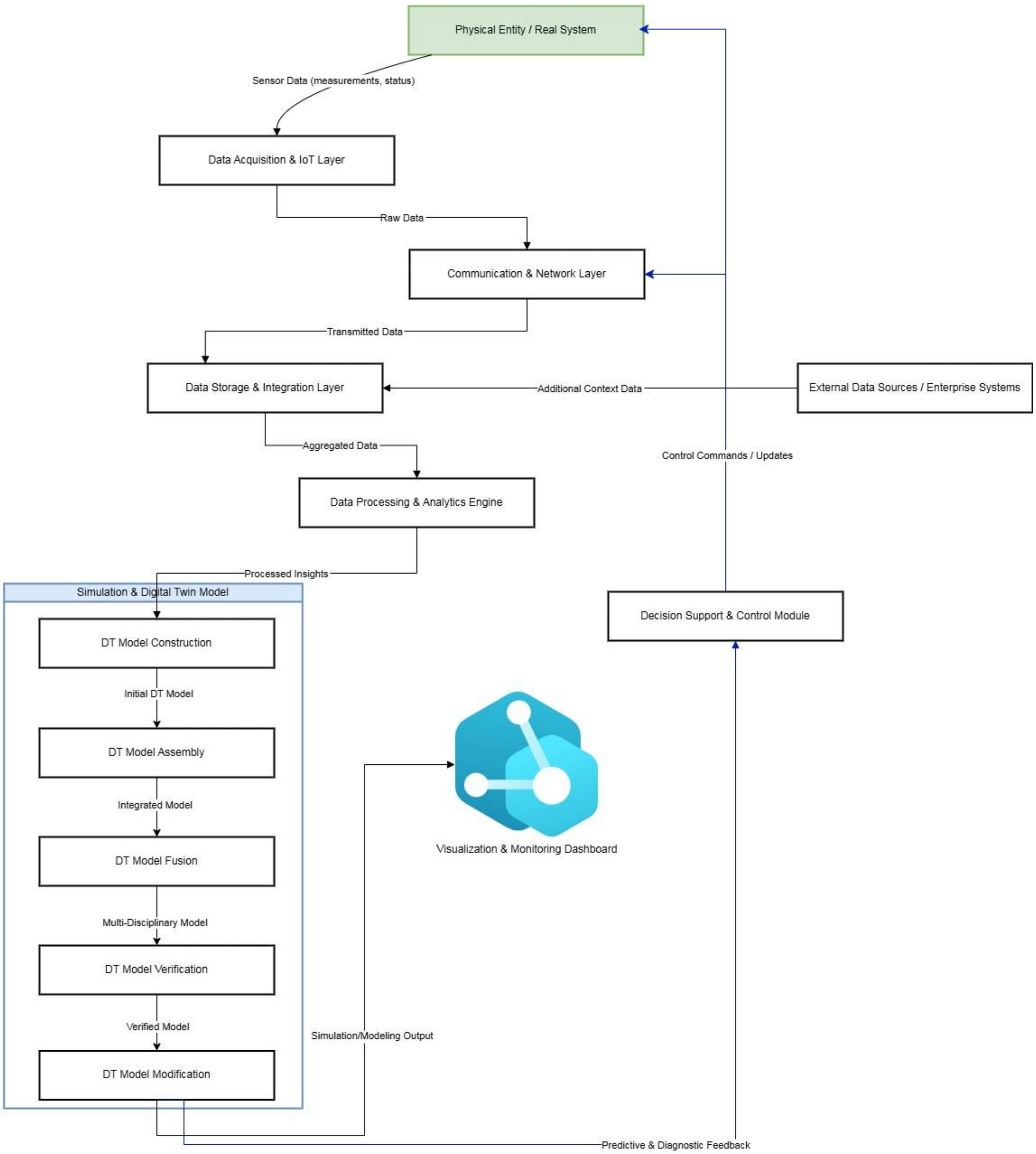

*Figure 6. General DT Architecture*

# 5. Discrete Event Simulation (DES)

In this section, we overview the strengths and weaknesses of Discrete Event Simulation (DES) as a modeling method for representing Service and Manufacturing Systems at the Product/Production Line, Production Facility, and Enterprise Levels. DES is a modeling technique widely used to study the behavior of physical systems where events occur at discrete points in time, and the system's state changes only at these events [92]. It is particularly effective for systems with queues, such as production lines in manufacturing or service centers like call centers, where the timing and sequence of events are critical.

At the Product/Production Line Level, DES is commonly applied to model individual products or specific production lines, focusing on detailed operations such as machine processing times, setup times, and queue behaviors. This level is critical for optimizing production efficiency and identifying bottlenecks [93]. Table 3 summarizes the strengths and weaknesses of DES at the product/production line, production facility, and enterprise levels.

*Table 3. Strengths and weaknesses of DES at different levels*

| **Level** | **Strengths** | **Weaknesses** |
|---|---|---|
| Product/Production Line | ➢ Simulation and analysis of the production line with different scenarios [94]<br>➢ Balancing the production line [95]<br>➢ Modeling the production line [96, pp. 385–390], [97]<br>➢ Bottleneck detection [93], [95], [98], [99]<br>➢ Evaluation of equipment utilization rate and failure of a production line [99]<br>➢ Verification of the line capacity and determining optimal buffer placement and size [98]<br>➢ Evaluating "What-if" Scenarios [100] | ➢ Dependency on the user for interpretation of quantitative results to get qualitative outputs [94], [98]<br>➢ Sensitive to input data [94], [95], [98]<br>➢ The level of detail required can be hard to define [98] |
| Production Facility | ➢ Simulation of the plant [99]<br>➢ Simulation of material handling systems within a factory [95]<br>➢ Analyzing the effects of plant layout on production efficiency [95]<br>➢ Detecting bottlenecks not just at a single line but across the entire facility [95], [101]<br>➢ Evaluation of logistic planning and design within the facility [96, pp. 391–396]<br>➢ Studying the shifts and labor movement across the facility [95]<br>➢ Analyzing the inventory management policies on the overall facility performance [95]<br>➢ Optimization of facility layouts [101] | ➢ Labor-intensive and time-consuming model development [95], [98], [102]<br>➢ User-dependent for simulation and analysis [98]<br>➢ Simulation is Not Optimization. Finding the optimal solution requires coupling DES with optimization algorithms, which adds complexity [95], [98]. |
| Enterprise | ➢ Can be used as a decision support tool at the enterprise level [95]<br>➢ Analyzing inter-enterprise integration, such as between suppliers and manufacturers [95]<br>➢ Can help evaluate how well certain strategic objectives are met, such as manufacturing lead time [95] | ➢ Highest complexity in the integration of data and processes (with other enterprise systems (e.g., ERP, MES)) [95], [98]<br>➢ Challenges related to gathering comprehensive and accurate data across the entire enterprise [95], [96, pp. 219–230], [98] |

# 6. Real-World Implementations

The practical application of digital twins at different levels has shown significant benefits in improving operational visibility, predicting maintenance needs, and making strategic decisions. However, challenges unique to each level exist. These include high initial costs and data integration problems at the production line level, increased model complexity and cybersecurity risks at the facility level, and oversimplification and organizational obstacles at the enterprise level.

| Level | Company- Application | Strengths (in general) | Weaknesses (in general) |
|---|---|---|---|
| Product/Production Line | ➢ Rolls-Royce (UK) - Jet Engines [72], [103]<br>➢ STEP Tools Inc. - Digital Twin Machining [104]<br>➢ GE- Jet Engines [105]<br>➢ Tesla- Car [106] | ➢ Identification of fault parameters and performing quantitative diagnosis [107]<br>➢ Allows for offline programming, potentially reducing downtime [74]<br>➢ Enables real-time quality inspection [107]<br>➢ Reducing reactive maintenance by 40% (for GE Engines) [105] | ➢ Variability in manufacturing vendor products can lead to data heterogeneity, posing a challenge to integration [73]<br>➢ Synchronization and consistency [74]<br>➢ Fidelity of models [74]<br>➢ Communication latency between physical and virtual components [74] |
| Production Facility | ➢ GE- Wind Farm [108]<br>➢ Volkswagen- An entire plant [30], [109]<br>➢ SIEMENS (Germany)- Programmable Logic Controllers (PLCs) [105]<br>➢ BMW- Car Factory [110], [111] | ➢ Reducing production [109]<br>➢ Quality improvement [109]<br>➢ Overall Equipment Effectiveness (OEE) improvement [74]<br>➢ Self-Organizing factory environment [107]<br>➢ Analyzing production bottlenecks [107] | ➢ Complexity of holistic models [74]<br>➢ Cost of implementation [74]<br>➢ Data governance and management at the factory level [28] |
| Enterprise | ➢ First Abu Dhabi Bank- IT [112] (using SAP LeanIX [113]) | ➢ Supply Chain integration [73]<br>➢ Cross-system and cross-platform interoperability [73] | ➢ Lack of standardized connection and communication [74]<br>➢ High expense and time [114] |

## 7. Roadblocks to Full Potential

Although many researchers have considered developing DT frameworks and models in the laboratories, there are still several roadblocks that inhibit Digital Twins from reaching their full potential. In this section, some of the most important of those barriers are addressed.

- **Data-Related Issues**: Variability in manufacturing vendor products can lead to *data heterogeneity*, posing a challenge for integration within a DT [73]. Also, a lack of device communication and data collection standards can compromise the *quality of data* being processed for DT, affecting its performance [28]. Furthermore, *managing the large volumes of data* involved in DTs, including identification, access, transformation, and quality assurance, can be a significant challenge [28]. In addition, ensuring data security and integrity is still a roadblock to DT implementation [28], [114], [115, p. 23].

- **Expensive Investments and Need for Cost-Benefit Analysis**: Implementing DT can involve expensive investments. Industries need to perform thorough cost-benefit analyses before committing to full-scale DT implementation, especially considering the potentially high costs and long timelines for complex systems. [28]. Many companies have expectations of a quick return on investment that might be unrealistic [74].

- **Lack of Standards and Regulations**: Widespread adoption is impeded by the lack of an industry-wide framework for DT, which includes standards for uniformity and a common understanding of interfaces. Despite ongoing efforts to standardize, such as ISO 23247 for manufacturing, these initiatives are still in the early stages of development. [28].

- **Challenges in Achieving High-Fidelity Models and Real-Time Synchronization**: A DT needs to be a near-identical copy of its physical counterpart with high accuracy and near real-time synchronization [28]. A major technical problem is achieving and sustaining this high fidelity, particularly for complex systems impacted by unpredictability and uncertainties [74].

- **Overlooking the Human Element**: The majority of DT studies have concentrated on industrial assets, ignoring the value of human interaction and the possibility of "Digital Twins for people" in the production setting [74], [107]. Further study is needed to integrate human elements and create human-centered human-machine collaboration tactics through DTs [107].
- **Technological maturity**: While technologies like IoT, AI, and simulation are evolving, their seamless integration and specific tools for DT development are still maturing [28], [114].
- **Limited Practical Validation**: There is a lack of extensive real-world validations of the proposed DT implementation models in actual manufacturing processes [74].